\documentclass[12pt]{article}
\usepackage{amsmath,epsfig,subfigure}
\usepackage{amsfonts,amsopn}
\usepackage{cases}

\newtheorem{thm}{Theorem}[section]

\newtheorem{lem}[thm]{Lemma}

\newtheorem{rem}[thm]{Remark}
\newtheorem{defn}[thm]{Definition}

\newcommand{\thmref}[1]{Theorem~{\rm \ref{#1}}}
\newcommand{\lemref}[1]{Lemma~{\rm \ref{#1}}}

\makeatletter \@addtoreset{equation}{section}

\newcommand{\beq}[1]{\begin{equation} \label{#1}}
\newcommand{\eeq}{\end{equation}}
\newcommand{\bed}{\begin{displaymath}}
\newcommand{\eed}{\end{displaymath}}

\newcommand{\bedd}{\bed\begin{array}{l}}
\newcommand{\eedd}{\end{array}\eed}
\newcommand{\bea}{\bed\begin{array}{rl}}
\newcommand{\eea}{\end{array}\eed}

\newcommand{\disp}{\displaystyle}
\newcommand{\ad}{&\!\!\!\disp}
\newcommand{\aad}{&\disp}
\newcommand{\barray}{\begin{array}{ll}}
\newcommand{\earray}{\end{array}}

\newcommand{\lbar}{\overline}

\newcommand{\M}{{\cal M}}
\newcommand{\be}{\beta}
\newcommand{\dl}{\delta}
\newcommand{\e}{\varepsilon}
\newcommand{\al}{\alpha}
\newcommand{\De}{\Delta}
\newcommand{\ga}{\gamma}

\newcommand{\sg}{\sigma}

\newcommand{\rr}{{\Bbb R}}
\newcommand{\nn}{\Bbb N}
\newcommand{\zz}{\mathbb{Z}}
\newcommand{\one}{{\rm 1}\hspace*{-0.035in} {\rm l}}
\newcommand{\cd}{(\cdot)}
\newcommand{\ar}{\rightarrow}
\newcommand{\F}{{\cal F}}
\newcommand{\Ll}{{\cal L}}

\newcommand{\p}{\mathfrak{p}}
\newcommand{\m}{\mathfrak{m}}

\newcommand{\nd}{\noindent}

\newcommand{\diag}{{\rm{diag}}}
\def\para#1{\vskip 0.1\baselineskip\noindent{\bf #1}}
\def\qed{\strut\hfill $\Box$}

\def\({\left(}
\def\){\right)}

\newcommand{\wdt}{\widetilde}
\newcommand{\sui}{\sum_{\iota \in \ZZ}}
\newcommand{\oal}{{I_{\{\al =\iota\}}}}
\newcommand{\ZZ}{{\mathbb Z_+}}

\begin{document}
\title{Hybrid Competitive Lotka-Volterra Ecosystems: Countable Switching States
and Two-time-scale Models\thanks{This research was supported in part by the National Science Foundation under grant DMS-1207667.}}
\author{Trang Bui\thanks{Department of Mathematics, Wayne State University, Detroit, MI 48202. Email: trang.bui@wayne.edu.} \and G. Yin\thanks{Department of Mathematics, Wayne State University, Detroit, MI 48202.
Email: gyin@math.wayne.edu.} }

\maketitle
\begin{abstract}
This work is concerned with  competitive Lotka-Volterra model
with Markov switching.
A novelty of the contribution is that the Markov chain has a countable state space.
Our main objective of the paper is
 to reduce the computational complexity by using the two-time-scale systems.
Because
existence and uniqueness  as well as continuity of  solutions for
 Lotka-Volterra ecosystems with Markovian switching in which the switching takes place in a countable set are not available, such properties are  studied first.
The two-time scale feature is highlighted by introducing a small parameter into the generator of the Markov chain. When the small parameter goes to 0, there is a limit system or reduced system.
It is established in this paper that if the reduced system possesses certain properties such as permanence and extinction, etc., then the  complex system also has the same properties when the parameter is sufficiently small.
These results are obtained by using the perturbed Lyapunov function methods.

\vskip 0.25 true in
\nd{\bf Key Words.} Lotka-Volterra model, perturbed Lyapunov method, Markov chain, singular perturbation, extinction, permanence, switching diffusion.

\vskip 0.25 true in
\nd {\bf Mathematics Subject Classification.}  60J60, 60J27, 92D25, 92D40.
\end{abstract}

\newpage

\setlength{\baselineskip}{0.22in}
\section{Introduction}
 Introduced by Lotka \cite{Lotka} and Volterra \cite{Vol}, the well-known Lotka-Volterra models have been  investigated extensively in the literature
  and used widely in ecological and population dynamics, among others. When two or more species live in close proximity and share the same basic requirements, they usually compete for resources, food, habitat, or territory.
Initially being posed as a deterministic model, subsequent study has taken randomness into consideration; see
\cite{DDNY,HNY,NY17} and references therein.
Recent effort on the so-called hybrid systems has much enlarged the applicability of Lotka-Volterra systems. One class of such hybrid systems uses a continuous-time Markov chain to model environmental changes and other random factors not represented in the usual stochastic differential equations; see \cite{YinZ10}
for a comprehensive study of switching diffusions.

Deterministic Lotka-Volterra systems have been studied by many people. A number of important results were obtained.
 A set of sufficient conditions for the existence of a globally stable equilibrium point in various
 models of the $n$-dimensional Lotka-Volterra system was obtained in \cite{Li-Tang}; limit cycles for some deterministic competitive three-dimensional Lotka-Volterra systems were treated in \cite{Xiao}. Random perturbations to the Lotka-Volterra model were considered in the literature; see for example \cite{Arnold,Khas} and many references therein, and also \cite{Imhof, Khas2} for up-to-dated progress on stochastic replicator dynamics.
Recently, much effort has been devoted to studying the stochastic Lotka-Volterra with regime-switching; see  \cite{Du,LuoMao,KyY,YinZ, YinZ2}.
 While most  recent works focus on Markov chain with a finite state space,
 to take into consideration of various factors, it is also natural to consider the Markov chain with a countable state space, which is the effort of
the current paper.

One of our main aims here is to reduce the computational complexity.
Assuming that the Markov chain has a two-time-scale structure so that it has a fast changing part and a slowly varying part,
we show that the systems under consideration can be reduced to a
much simpler limit system.
In fact, the limit system is a diffusion that no switching is involved.
The rationale is that the original system with the small parameter is much more difficult to deal with, but the limit system is substantially simpler.
Using  properties such as permanence, stochastic boundedness, extinction etc.
 of the limit diffusion system, we can make inference about the more complex original system. The complexity reduction is achieved by using two-time-scale formulation and perturbed Lyapunov function methods.
 This line of thinking goes back to the work of \cite{BP78}, which has been much expanded to more general setting in \cite{Kush}; see also \cite{YinZh}.

 Mathematically, the time scale separation is obtained by introducing a small parameter $\e$. As $\e\to 0$, we obtain a limit system. We then show if the limit system has certain properties, then the complex original system also preserves the same property for sufficient small $\e$. Note that for the regime-switching Lotka-Volterra system with a countable state space for the switching, such properties as existence and uniqueness of solution,  and moment bounds etc. are not yet available. Thus, we first establish these properties.

  We model the random environments
  (e.g., different
    seasons, changes in nutrition and food resources, and other random factors)
   in the ecological system by a continuous-time Markov chain $\al(t)$ with a countable state space $\ZZ= \{1,2, \dots\}$ and a generator
   $Q = (q_{\al\be})$ satisfying
  $q_{\al\be} \geq 0$ for $\al\in \ZZ$ and $\be\not = \al$,
   and $\sum_{\be=1}^{\infty} q_{\al\be}=0$ for each $\al\in \ZZ$.
A stochastic Lotka-Volterra system in random environments can be described by the following stochastic differential equation
  (in the Stratonovich sense) with regime switching
$$dx_i(t)=x_i(t) \bigg\{ \bigg[ b_i(\al(t)) - \sum_{j=1}^n a_{ij}(\al(t))x_j(t) \bigg] dt + \sigma_i(\al(t)) \circ dw_i(t) \bigg\}, \  i=1,\dots,n,$$
where $w\cd  =(w_1\cd,\dots,w_n \cd)' $ is an $n$-dimensional standard Brownian motion, $b(\al) = (b_1(\al),\dots,b_n(\al))'$, $A(\al)=(a_{ij}(\al))$, and $\Sigma(\al)=\diag(\sigma_1(\al),\dots,\sigma_n(\al))$ with $\al \in \M$ represent
intrinsic growth rates, the community matrices, and noise intensities in different external environments, respectively.
It is well known that the above stochastic differential equation in the Stratonovich sense
is equivalent to the system in the It\^{o} sense
\beq{Ito}
dx_i(t)=x_i(t) \bigg\{ \bigg[ r_i(\al(t)) - \sum_{j=1}^n a_{ij}(\al(t))x_j(t) \bigg] dt + \sigma_i(\al(t)) dw_i(t) \bigg\}, \  i=1,\dots,n,
\eeq
where
 $r_i(\al) := b_i(\al)+\frac{1}{2}\sigma_i^2(\al) \ \hbox{ for each } \ i=1,2,\dots,n.$
 In ecology and biology, one prefers to
start the formulation of stochastic Lotka-Volterra systems using calculus in the
Stratonovich sense because each term has its clear ecological meaning.  However, for the analysis, the It\^o calculus should be used.
Assume throughout the paper that  the Markov chain $\al\cd$ and the Brownian motion $w\cd$ are  independent. Without loss of generality, we also assume that the initial conditions $x(0)$ and $\al(0)$ are non-random.

Note that
 $(x(t),\al(t))$ is a Markov process, whose
 generator $\Ll$ is given as follows (see \cite[Chapter 2]{YinZ10} and also \cite{YinZ,YinZ2} for a definition of the generator of a Markov process). For any $V: \rr^n \times \ZZ \mapsto \rr$ with $V(\cdot, \al)$ being twice continuously differentiable with respect to the variable $x$ for each $\al \in \ZZ$, we define
\beq{gen} \barray
\Ll V(x,\al):= &\disp \sum_{i=1}^n \frac{\partial}{\partial x_i} V(x,\al) x_i \bigg(r_i(\al)-\sum_{j=1}^n a_{ij}(\al)x_j \bigg) + \frac{1}{2} \sum_{i=1}^n \frac{\partial^2}{\partial x_i^2} V(x,\al) x_i^2 \sigma^2_i(\al)\\
& + \disp \sum_{\beta\in \ZZ, \beta \neq \al} q_{\al \beta} [V(x,\beta) - V(x,\al)].
\earray \eeq
Comparing to \cite{YinZ,YinZ2}, the Markov chain takes values in a countably infinite set.

In this paper, we focus on the study of hybrid Lotka-Volterra systems involving a
two-time-scale Markov chain
 with the help of the study of asymptotic properties of two-time-scale Markov chains  \cite{YinZh}.
In contrast to the existing results, our contributions are as follows. (i) We model the Lotka-Volterra ecosystems using hybrid systems in which continuous states (diffusion) and discrete events (switching) coexist and interact. A distinct feature of the modeling point is that the random discrete events take values in a countably infinite set. (ii)
Prior to this work,  existence and uniqueness of solution, continuity of sample paths, and stochastic boundedness of regime-switching Lotka-Volterra system with random switching taking values in a countable state space were not available. Our paper establishes these properties.
Although general regime-switching diffusions were considered in \cite{NY16}, the spatial variable $x$ there lives in the whole space $\rr^n$, whereas for the Lotka-Volterra systems considered here, $x\in \rr^n_+$. It needs to be established that the solution is in $\rr^n_+$ as well.
 (iii) This paper provides a substantial reduction of complexity. The two-time scale system is a system involving countably infinitely many equations, whereas the limit system is a single diffusion. Using the limit system as a bridge, we then obtain
 for example, if the limit system is stochastically bounded, or
  permanent, or going to extinction,  then the much more complex system
  with switching also preserves such properties as long as the parameter $\e>0$ is small enough.

We consider Lotka-Volterra ecosystems with a multiple number of species. The species competing against each other.
 Our main interests are to closely capture the dynamics, to reveal
whether or not the species  can exist, to find conditions that they will be permanent or  extinct. Such studies can be carried out using Lyapunov function
methods.
 The biological and ecological significance is as follows. For the rather complex ecosystems,
 we show that the complex models with countably infinite discrete states can be handled by using more manageable limit (reduced) diffusion system.  More specifically, the answers to such important questions as permanence and distinction of the competing species of the complex systems can be answered by examining the reduced systems. In doing so, we achieve a substantial reduction of computation complexity.

The rest of the paper is organized as follows. Section \ref{pros} studies existence, uniqueness, and continuity of solutions of the competitive Lotka-Volterra systems associated with a continuous-time Markov chains with a countable state space. We then introduce the Lotka-Volterra systems with two-time scales with the use of a singularly perturbed Markov chain and illustrate the properties or their solutions in Section \ref{pert}. We further provide the permanence and extinction of the systems with two-time-scale Markov chains through their limit systems in \ref{Extin} and \ref{Perm}. The paper is concluded with conclusions and remarks in Section \ref{Rem}.
Finally, an appendix containing the technical complements of
 the proofs of a number of technical results
  is provided.

\section{Existence, Uniqueness, and Continuity of Solutions}
 \label{pros}
Before getting to the two-time-scale systems, we first examine systems without the time scale separation. Existence, uniqueness, and continuity of solutions of the regime-switching Lotka-Volterra systems   when $\ZZ$ is countably infinite are not available. So we present these results first in what follows.
Denote
\beq{eq:notation}\barray \ad \Xi (x,\al) := \diag(x_1,\dots,x_n)[r(\al)-A(\al)x],\\
\ad \disp \xi_i(x,\al)= x_i ( r_i(\al) -\sum_{j=1}^n a_{ij}(\al) x_j ),\\
 \ad \disp s_i(x,\al) = x_i \sigma_i( \al)\\
 \ad \disp S(x,\al) =\diag(s_i(x,\al)).\earray\eeq

By a competitive system, we mean that all values in the community matrix $A(\al)$ are non-negative ($a_{ij}(\al) \geq 0$ for all $\al \in \ZZ$ and $i,j=1,2,\dots,n$). It is reasonable to assume that the competitions among the same species are strictly positive. Therefore, we assume
\begin{itemize}
\item[(A1)] For each $\al \in \ZZ = \{1,2,\dots\}$,
$a_{ii} (\al) >0$
and $a_{ij}(\al) \geq 0$ for $i,j =1,2,\dots,n$ and $j \neq i$.
\end{itemize}

In \cite{YinZ}, the existence and uniqueness for the switching diffusion model was obtained when the state space of the switching is finite. However, the state space of the Markov chain in our study is countable but not finite.
If $\al(t-) := \lim_{s \to t^-} \al(s) = \al$, then it can switch to $\beta$ at $t$ with intensity $q_{\al\beta}$. Denote for each $\al$,
 $q_\al = \sum_{\beta \in \ZZ,\beta \neq \al}q_{\al \beta}$.
 Note that $\al(t)$ may be written
  as the solution to a stochastic differential equation with respect to a Poisson random measure. To be more precisely, let $\p(dt,dz)$ be a Poisson random measure with intensity $dt \times \m(dz)$ and $\m$ be the Lebesgue measure on $\rr$ such that $\p(\cdot,\cdot)$ is independent of the Brownian motion $w(t)$.
  Using this fact, for each $\al \in \zz$, we can construct disjoint sets $\{ \Delta_{\al\beta}, \beta \neq  \al \}$ on the real line as follows
$$\barray
\De_{12} =[0,q_{12}),\\
\De_{13}=[q_{12},q_{12}+q_{13}),\\
\dots\\
\De_{21}=[q_1,q_1+q_{21}),\\
\De_{23}=[q_1+q_{21},q_1+q_{21}+q_{23}),\\
\dots
\earray$$
Define $h: \ZZ \times \rr \mapsto \rr$ by $h(\al,z)=\disp \sum_{\beta\in \ZZ,\beta \neq \al}^{\infty} (\be-\al) \one_{\{z \in \De_{\al\be}\}}$, where $\one_{\{z \in \De_{\al\be}\}}=1$ if $z \in \De_{\al\be}$ and $\one_{\{z \in \De_{\al\be}\}}=0$, is the indicator function. The process $\al(t)$ can be defined as a solution to
$$d\al(t) =\int_{\rr}h(\al(t-),z) \p(dt,dz),$$
where $\p(dt,dz)$ is a Poisson measure with intensity $dt \times \m(dz)$ and $\m$ is the Lebesgue measure on $\rr$.
We assume the following condition holds.

\begin{itemize}
\item[(A2)] The Markov chain having  generator $Q$ is strongly exponentially ergodic (see \cite{Ander})
in that
there exist a $K>0$ and a $\lambda_0>0$ such that
 \beq{exp-er}
 \sum_{\be=1}^\infty |p_{\al\be}(t)- \nu_\be |\le K \exp(-\lambda_0 t)\eeq
  for any positive integer $\al$ and $t>0$, where $\nu$ is the stationary distribution associated with the generator $Q$ and $p_{\al\be}(t) = P (\al(t) =\be| \al(0) =\al)$. Moreover,
\beq{A3}
M = \sup_{\al}{q_\al}=\sup_\al \sum_{\be\not =\al} q_{\al\be} < \infty.
\eeq
\end{itemize}

To proceed, we obtain the global solution in $\rr^n_+$ for the system.
Then we  establish the positivity of solution $x(t)$, finite moments, and continuity.  One of the main tools is to use an appropriate Lyapunov functions.
 The proofs of the results are relegated to the appendix for convenience.

\begin{thm}\label{supr}
Assume {\rm(A1)} and {\rm(A2)}. Then for any initial data
$x(0) = x_0 \in \rr_+^n$ and $\al(0)= \al \in \ZZ$,
there is a unique solution $x(t)=(x_1(t),\dots,x_n(t))'$ to \eqref{Ito} on $t \geq 0$, and the solution will remain in $\rr_+^n$ almost surely, i.e., $x(t) \in \rr_+^n$ a.s. for any $t \geq 0$.
\end{thm}

We next consider the stochastic boundedness. First, we recall the definition.

\begin{defn}{\rm
The solution $x(t)$ of \eqref{Ito} is stochastically bounded (or bounded in probability), if for any $\eta >0$, there is a constant $H=H_\eta$ such that for any $x_0 \in \rr^n_+$,
\beq{bounded}
\limsup\limits_{t \rightarrow \infty} P \{ \vert x(t) \vert \leq H\} \leq 1-\eta.
\eeq
}\end{defn}

\begin{thm}\label{moment} Under the conditions of  \thmref{supr} and for any $p>0$ satisfying
\beq{ment}
\disp \sup_{\al \in \zz^+} \sum_{i=1}^n  \dfrac{1+ pb_i (\al) +\frac{p^2}{2}\sigma_i^2 (\al)}{a_{ii} (\al)}  < \infty,
\eeq
we have
\begin{equation}
\sup_{t \geq 0} E\bigg[ \sum_{i=1}^n x_i^p(t) \bigg]  \leq K < \infty.
\end{equation}
\end{thm}

By virtue of Tchebychev's inequality, a direct consequence of \thmref{moment} is that the solution $x(t)$ is stochastically bounded.
Next we obtain the sample path continuity.

\begin{thm}\label{continu} The solution $x(t)$ to \eqref{Ito} is continuous a.s.
\end{thm}

\begin{rem}{\rm In fact, 
for almost all sample paths of the solutions \eqref{Ito} are H\"older continuous with exponent $\gamma < \frac{1}{4}$. That is, except a null set $N$ with probability 0, for all $\omega \in \Omega \backslash N$, there exists a random variable $h(\omega) >0$ satisfying \begin{equation}
P \Bigg\{ \omega: \sup \limits_{0 \leq s,t < \infty, |t-s| < h(\omega)} \frac{|x(t, \omega) - x(s, \omega)|}{|t-s|^\gamma} \leq \frac{2}{1-2^{-\gamma}} \Bigg\} =1;
\end{equation}
see the proof in the appendix for more details.
}\end{rem}

\section{Two-Time-Scale Models}\label{pert}
\subsection{Two-Time-Scale Markov Chains}
 Recall that a generator $Q$ or its corresponding Markov chain is said to be irreducible if the system of equations
\beq{irre}
\nu Q = 0,\ \ \ \ \sum_{\al=1}^{\infty} \nu_\al= 1
\eeq
has a unique solution $\nu = (\nu_1,\nu_2,\dots)'$ satisfying that $\nu_\al > 0$ for $\al=1,2,\dots$ Such a solution is termed a stationary distribution. Throughout the rest of the  paper, we assume that the Markov chain
 has a fast varying part and slowly varying part in that
$\al(t)=\al^\e (t)$, with generator
\beq{t-time-Q} Q^\e = \frac{Q}{\e} + Q_0,\eeq
where $Q$ is a generator of a Markov chain that is irreducible and $Q_0$ is a generator of another Markov chain. We do not have any restrictions on $Q_0$. For simplicity, we use \eqref{t-time-Q} in this paper. Although $Q_0$ appears in \eqref{t-time-Q}, the asymptotic properties are dominated by $Q$. It is possible to consider more complex models with more structure on $Q$; see also the concluding remark section for more discussion.
For the subsequent study, we need a couple of preliminary results. The proofs of (i) and (ii) in \lemref{add} can be found in \cite[Theorem 4.5, Theorem 4.48, Lemma 5.1]{YinZh}.
Denote $p^\e(t) = (P(\al^\e(t) = \al): \al=1,2,\dots)$, $p^\e_{\al\be}(t,t_0) = P (\al^\e(t) =\be | \al^\e(t_0) =\al)$, and $P^\e(t,t_0)$ is the transition
matrix $(p^\e_{\al\be}(t,t_0))$.

\begin{lem}\label{add}
Assume that  for $Q$ given in \eqref{t-time-Q} satisfies
{\rm (A2)}.
Then
there exists a positive constant $\kappa_0$ such that
\begin{enumerate}
\item[{\rm(i)}] For the probability distribution vector $p^\e(t) \in \rr^{1 \times \infty}$
\begin{equation}
p^\e(t)=  \nu + O(\e  + e^{-\kappa_0 t /\e})
\end{equation}
uniformly in $(0,t)$.
\item[{\rm(ii)}] For the transition probability matrix $P^\e(t,t_0)$, we have
\begin{equation}
P^\e(t,t_0) = P_0(t) +O\Big( \e + e^{-\kappa_0 (t-t_0)/\e}\Big),
\end{equation}
uniformly in $(t_0,t)$, where $P_0(t) = \one \nu$.
with $\one =(1,1,\dots)'$ being an infinite column vector having
all entries $1$, and $\nu=(\nu_1,\nu_2,\dots)$ is the row vector of stationary distribution associated with the Markov chain with generator $Q$.
\end{enumerate}
\end{lem}

\begin{thm}\label{O-e} Assume \eqref{irre}. Then for each $\al=1,2,\dots$,
\beq{meas}
E \Bigg[\int_0^\infty e^{-t} (I_{\{\al^\e(t) =\al\}}-\nu_\al) dt\Bigg]^2 = O(\e).
\eeq
\end{thm}

Let $x^\e(t) \in \rr^n$ for $t \geq 0$ be given by
\beq{perturb}
dx^\e(t)=\diag(x^\e_1(t),\dots,x^\e_n(t))\bigg[\left(r(\al^\e(t))-A(\al^\e(t))x^\e(t)\right)dt+
\Sigma(\al^\e(t)) dw(t) \bigg],
\eeq
 with the initial conditions $x(0) = x_0$ and $\al^\e(0) = \al_0 \in \ZZ$.
Under  (A1) and (A2), we can construct the solutions of the two-time-scale stochastic differential equations by using similar method as in \thmref{supr}. The existence and uniqueness of solutions of the stochastic differential equations \eqref{perturb} hold; $0 \in \rr^n$ is a stationary point for each equation in \eqref{perturb}.

\begin{rem}\label{about-sec2}{\rm
Note that existence and uniqueness of solutions, path continuity, and moment bounds established in Section \ref{pros} hold for the two-time scale system \eqref{perturb}. Our main effort below
is to show how we may reduce the computational complexity.
}
\end{rem}

\begin{lem} Under {\rm(A1)} and {\rm(A2)},  $\{x^\e \cd\}$  given by \eqref{perturb} converges weakly to $x \cd$ such that $x\cd$ satisfies
\beq{aver}
dx(t)= \overline{\Xi}(x(t))dt+ \overline{\Lambda}(x(t))dw(t),
\eeq
where $\overline{\Xi}(x)= \sum_{\al=1}^{\infty} \Xi(x, \al) \nu_\al $, $\overline{\Lambda}(x) \overline{\Lambda}'(x)= \sum_{\al=1}^{\infty} S(x,\al) S'(x,\al) \nu_\al $.
\end{lem}

Weak convergence of $x^\e\cd$ to $x\cd$ is a basic notion in stochastic processes. A definition can be found in \cite[pp.371-376]{YinZ10}.
For convenience, we denote
\beq{bar-sys}\barray \ad \overline{r}_i= \disp \sum_{\al=1}^\infty r_i(\al) \nu_\al, \  \lbar{b}_i =\disp \sum_{\al=1}^\infty b_i(\al) \nu_\al,\
 \overline{a}_{ij} = \disp \sum_{\al=1}^\infty a_{ij} (\al) \nu_\al,\
 \overline{\sigma}_i =\disp\sqrt{ \sum_{\al=1}^\infty {\sigma^2_i(\al) \nu_\al}} \ \hbox{ and}\\
 \ad \lbar{\xi}_i (x) =\sum_{\al=1}^\infty \xi_i(x,\al) \nu_\al,\  \disp \lbar{\lambda}_i(x) = x_i \lbar{\sigma}_i = \sqrt{ \sum_{\al=1}^\infty s_i^2(x,\al) \nu_\al} ,\earray\eeq
 where $s_i(x,\al), S(x,\al)$ and $\Xi(x,\al)$ are defined in \eqref{eq:notation}.
The averaged system can be written
component-wise
as
\beq{c-aver}
dx_i(t) = x_i(t) \Bigg\{ \bigg[\lbar{r}_i -\sum_{j=1}^n \lbar{a}_{ij} x_j(t)\bigg]dt + \lbar{\sigma}_i dw_i(t)\Bigg\}.
\eeq

\begin{rem}{\rm Note the following facts.
\begin{itemize}
\item The proof of the above lemma is similar to the development in \cite[Ch.8]{YinZh}.
\item The averaged system \eqref{aver} is a Lotka-Volterra diffusion system, whose coefficients are an average with respect to the stationary measure $\nu$.
    Hence, under (A1) and (A2),  we can prove that the averaged system \eqref{aver} has a unique solution that is continuous together with moment bounds. This follows the way of treating nonlinear stochastic differential equations. First, we show that there is a local solution and then extend the solution to a global solution by using stopping time argument; see for example, \cite[Theorem 2.1]{MMR}.

\end{itemize}
}\end{rem}

In the study of stochastic population systems, we are interested in the permanence and extinction of the population.
We shall study this by means of the corresponding limit system.
Treating directly stability of dynamic systems containing two-time-scale
 Markov chains is a complex matter.
 However, considering this problem using limit system is much simpler.
 Some
  earlier work concerning the stability of those systems can be found in \cite{Bado}.
In this study, our goal is to establish the permanence and extinction of \eqref{perturb} for sufficiently small $\e$.
Here, from a Lyapunov function $V(x)$ of the averaged system, we construct a perturbed Lyapunov function for the
 more complex original
 system containing the fast varying Markov chain.
 The method we use is motivated by arguments in  \cite[pp. 148-149]{Kush}. The averaged system is a diffusion without switching, whereas in the original system, the switching states belong to a countably infinite set.
 Using the limit system, we can examine the original system, which is much easier that dealing with the original system directly.
 As a result, our approach leads to a significant reduction of complexity.

\subsection{Preliminary Calculations}
To proceed, we first present some preliminary calculations using perturbed Lyapunov function for preparation on study of various properties of the complex original system.
Let $\F^\e_t = \sigma \{ x^\e(s), \al^\e (s), s \leq t \}$, and $E^\e_t$ be the expectation conditioned on $\F^\e_t$. For a suitable function $\zeta (t)$, define the operator $\Ll ^\e$ by
\beq{def}
\Ll^\e \zeta(t) = \lim_{\delta \downarrow 0} \frac{1}{\delta} E^\e_t[\zeta(t+\delta) -\zeta(t)].
\eeq
The generator is as defined in \eqref{gen} with the switching part given  by \eqref{t-time-Q}.
As a result,  the generator of the switching diffusion
process is $\e$ dependent.
Let $V(x)$ be a Lyapunov function associated with the averaged system \eqref{aver}
independent of the discrete component.
Using \eqref{gen} for $V(x^\e(t))$ where $x^\e(t)$ is the solution of system \eqref{perturb}, we obtain
$$\Ll^\e V(x^\e(t)) = \sum_{i=1}^n V_{x_i}(x^\e(t)) \xi_i(x^\e(t), \al^\e(t)) + \frac{1}{2}\sum_{i=1}^n V_{x_i x_i} (x^\e(t)) s_i^2(x^\e(t), \al^\e(t)).$$
Define
\beq{per-term}
V_1^\e(x,  t) = E^\e_t \int_t^\infty e^{t-u} \sum_{i=1}^n V_{x_i} (x) [\xi_i(x, \al^\e(u)) - \lbar{\xi}_i(x)]du.
\eeq
\beq{per2}
V^\e_2(x,  t) = E^\e_t \int_t^\infty e^{t-u} \frac{1}{2}\sum_{i=1}^n V_{x_i x_i}(x )[s_i^2(x, \al^\e(u)) - \lbar{\lambda}_i^2(x)]du.
\eeq
This implies that with $\al^\e(t)= \ell$,
\bea V^\e_1(x,  t)\ad = E^\e_t \int_t^\infty e^{t-u} \sum_{i=1}^n V_{x_i }(x )
\sum_{k=1}^\infty \xi_i(x, k) [I_{\{\al^\e(u) = k\}} - \nu_k]du \\
\ad =\int^\infty_t e^{t-u}\sum_{i=1}^n V_{x_i}(x)\sum^\infty_{k=1}\xi_i(x, k)[p_{\ell k}(u)- \nu_k] du .\eea
Hence,
\beq{e-V1}V_1^\e (x^\e(t),  t) = O(\e) [V(x^\e(t))+1].
\eeq
To proceed, we use a notation
\beq{e-G} \barray
G(x,\al)
&=\disp \sum_{i=1}^n \sum_{k=1}^\infty V_{x_i}(x) \xi_i(x, k) \bigg[ I_{\{\al =k\}} - \nu_k \bigg].
\earray \eeq
Using \eqref{def}, \eqref{e-V1}, and \eqref{e-G},
\beq{egen}\barray
&\!\!\!\Ll^\e V_1^\e(x^\e(t), t)\\
 &= \lim \limits_{\delta \downarrow 0} \frac{1}{\delta} E^\e_t [V_1^\e(x^\e(t+ \delta), t+\delta)- V_1^\e(x^\e(t),  t)]\\
 &= \disp \sum_{k=1}^\infty \sum_{i=1}^n \sum_{j=1}^n [V_{x_i}(x^\e(t)) \xi_i(x^\e(t),k)]_{x_j} \xi_j(x^\e(t), \al^\e(t)) E^\e_t \int_t^\infty e^{t-u} [ I_{\{\al^\e(u) =k\}} - \nu_k] du\\
 &\ \disp + O(\e) (V(x^\e(t))+1) - \sum_{i=1}^n V_{x_i}(x^\e(t))[ \xi_i(x^\e(t), \al^\e(t))- \lbar{\xi}_i(x^\e(t))]\\
 &\disp = O(\e) (V(x^\e(t))+1) - \sum_{i=1}^n V_{x_i}(x^\e(t)) [\xi_i(x^\e(t), \al^\e(t))-\lbar{\xi}_i(x^\e(t))].
\earray \eeq
Similar to the estimate of $V^\e_1 (x,  t)$, it can be verified that
\beq{e-V2}\barray
& V^\e_2(x^\e(t),  t) = O(\e) [V(x^\e(t))+1],\\
& \disp \Ll^\e V^\e_2(x^\e(t),  t) = O(\e) [V(x^\e(t))+1]
 - \frac{1}{2}\sum_{i=1}^n V_{x_i x_i} (x^\e(t)) [s_i^2(x^\e(t), \al^\e(t))-\lbar{\lambda_i}^2(x^\e(t))].
\earray \eeq
Define $V^\e(x^\e(t),t) = V(x^\e(t)) + V^\e_1(x^\e(t), t) + V^\e_2(x^\e(t),  t)$ satisfying the following properties:
\bea \ad V^\e(x^\e(t),t) = V(x^\e(t)) +O(\e) (V(x^\e(t))+1).\\
\ad \Ll^\e V^\e(x^\e(t), \al^\e(t), t) = O(\e) (V(x^\e(t))+1) + \Ll V(x^\e(t), \lbar{\al}),\eea
 where $$\Ll V(x^\e(t), \lbar{\al})= \disp \sum_{i=1}^n V_{x_i}(x^\e(t)) x_i^\e(t) \bigg( \lbar{r}_i - \sum_{j=1}^n \lbar{a}_{ij} x^\e_j(t)\bigg) + \frac{1}{2} \sum_{i=1}^n V_{x_i x_i}(x^\e(t)) [x_i^\e(t)]^2 \lbar{\sigma}^2_i.$$

\subsection{Stochastic Boundedness}\label{stb}

First, under suitable conditions,
the averaged system is stochastically bounded.
This follows from a specialization of the
proof of \cite[Theorem 3.1]{YinZ} (for the case that
the switching set has only one element), which is a refinement of the
 arguments of moment bounds in \cite{Mao}.

\begin{lem}\label{bnd} Assume that {\rm(A1)}, {\rm(A2)}, and \eqref{irre} are satisfied. Then the following statements hold for the solution $x(t)$ of \eqref{aver}.
\begin{itemize}
\item[{\rm(1)}] For any $p>0$,
\begin{equation}
\sup_{t \geq 0} E\bigg[ \sum_{i=1}^n x_i^p(t) \bigg]  \leq K < \infty.
\end{equation}
\item[{\rm(2)}] For any $p>0$,
\begin{equation}
\limsup_{t \to \infty} E \bigg[|x(t)|^p \bigg] \leq K < \infty.
\end{equation}
\item[{\rm(3)}]  The solution of the averaged system \eqref{aver}, namely, $x(t)$, is stochastically bounded, i.e., for any $\dl >0$, there is a constant $H = H(\dl) $ such that for any $x_0 \in \rr ^n_+$, we have
\begin{equation}
\limsup_{t \to \infty} P \{ |x(t)| \leq H \} \geq 1-\dl.
\end{equation}
\end{itemize}
\end{lem}

With the lemma above, we proceed to show that the solution of system \eqref{Ito} also has the same boundedness property
if $\e$ is small enough. Note that the next theorem should be compared with \thmref{moment}. Different from \thmref{moment}, the condition \eqref{ment} is not needed in the following theorem. However, it is required that $\e$ be small enough.

\begin{thm}
\label{Bound} Assume that {\rm(A1)}, {\rm(A2)}, and \eqref{irre} are satisfied. Then the following statements hold for the solution $x^\e(t)$ of \eqref{perturb} for $\e$ sufficiently small.
\begin{itemize}
\item[{\rm(1)}] For any $p>0$,
\begin{equation}
\sup_{t \geq 0} E\bigg[ \sum_{i=1}^n [x_i^\e(t)]^p \bigg]  \leq K < \infty.
\end{equation}
\item[{\rm(2)}]  For any $p>0$,
\begin{equation}
\limsup_{t \to \infty} E \bigg[|x^\e(t)|^p \bigg] \leq K < \infty.
\end{equation}
\item[{\rm(3)}] The process $x^\e(t)$ is stochastically bounded. That is, for any $\dl >0$, there is a constant $H = H(\e, \dl)$ such that for any $x^\e_0 \in \rr ^n_+$, we have
\begin{equation}
\limsup_{t \to \infty} P \{ |x^\e(t)| \leq H \} \geq 1-\dl.
\end{equation}
\end{itemize}
\end{thm}

\para{Proof.} We use perturbed Lyapunov function methods to prove this theorem.
Consider \bea
\ad \disp V(x)=\sum_{i=1}^n [x_i]^p \ \hbox{ and }\ \
 \disp \tilde{V}(x) = \sum_{i=1}^n \bigg[ [x_i]^\ga -1 -\ga \log x_i \bigg].\eea
Similar to \eqref{per-term} and \eqref{per2}, we define
\bea
& V^\e_1(x,t) =E^\e_t \disp \int_t^\infty e^{t-u} \disp\sum_{i=1}^n V_{x_i} (x) [\xi_i (x, \al^\e(u))- \lbar{\xi}_i(x)]du,\\\\
&V^\e_2(x,t) = E^\e_t \disp\int_t^\infty e^{t-u} \dfrac{1}{2} \disp\sum_{i=1}^n V_{x_i x_i} (x) [s^2_i (x, \al^\e(u)) -\lbar{\lambda}_i^2(x)]du,\\\\
& \tilde{V}^\e_1(x,t) =E^\e_t \disp \int_t^\infty e^{t-u} \disp\sum_{i=1}^n \tilde{V}_{x_i} (x) [\xi_i (x, \al^\e(u))- \lbar{\xi}_i(x)]du,\\\\
&\tilde{V}^\e_2(x,t) = E^\e_t \disp\int_t^\infty e^{t-u} \dfrac{1}{2} \disp\sum_{i=1}^n \tilde{V}{x_i x_i} (x) [s^2_i (x, \al^\e(u)) -\lbar{\lambda}_i^2(x)]du,
\eea
which have the following properties:
\bea
 \ad V^\e_1(x^\e(t),t) = O(\e) [V(x^\e(t)) +1], \\ \ad \tilde{V}^\e_1(x^\e(t), t) = O(\e) [\tilde{V}(x^\e(t)) +1], \\
\ad V^\e_2(x^\e(t), t) = O(\e) [V(x^\e(t)) +1], \\
\ad \tilde{V}^\e_2(x^\e(t),t) = O(\e) [\tilde{V}(x^\e(t)) +1] ,\\
\ad \Ll^\e V^\e_1(x^\e(t), t) = O(\e)(V(x^\e(t))+1)- \disp \sum_{i=1}^n V_{x_i} (x^\e(t))[\xi_i(x^\e(t), \al^\e(t)) -\lbar{\xi}_i (x^\e(t))],\\
\ad
\Ll^\e \tilde{V}^\e_1(x^\e(t), t) = O(\e)(\tilde{V}(x^\e(t))+1)- \disp \sum_{i=1}^n \tilde{V}_{x_i} (x^\e(t))[\xi_i(x^\e(t),\al^\e(t)) -\lbar{\xi}_i (x^\e(t))],\\
\ad \Ll^\e V^\e_2(x^\e(t),  t) = O(\e)(V(x^\e(t))+1)- \disp \sum_{i=1}^n V_{x_i x_i} (x^\e(t))[s^2_i(x^\e(t),\al^\e(t)) -\lbar{\lambda}^2_i (x^\e(t))],\\
\ad \Ll^\e \tilde{V}^\e_2(x^\e(t), t) = O(\e)(\tilde{V}(x^\e(t))+1)- \disp \sum_{i=1}^n \tilde{V}_{x_i x_i} (x^\e(t))[s^2_i(x^\e(t),\al^\e(t)) -\lbar{\lambda}^2_i (x^\e(t))].
\eea
Define
\bea \ad V^\e (x^\e(t), t) = V(x^\e(t)) + V^\e_1 (x^\e(t),t)+V^\e_2 (x^\e(t), t)\ \hbox{
 and }\\
\ad \tilde{V}^\e (x^\e(t),t) = \tilde{V}(x^\e(t)) + \tilde{V}^\e_1 (x^\e(t),t)+\tilde{V}^\e_2 (x^\e(t), t).\eea Then
\beq{tV-e}
\barray
\Ll^\e \tilde{V}^\e(x^\e(t), t)& \leq O(\e) (\tilde{V}(x^\e(t))+1) + \ga \disp \sum_{i=1}^n \bigg\{ -\lbar{a}_{ii} [x^\e_i(t)]^{\ga+1} + \bigg(\lbar{b}_i +\frac{\ga}{2} \lbar{\sigma}_i^2 \bigg)[x^\e_i(t)]^\ga\\
&\qquad + \disp \bigg(\sum_{j=1}^n \lbar{a}_{ji} \bigg) x^\e_i(t) -\lbar{b}_i \bigg\} \\
& \leq O(\e) \tilde{V}(x^\e(t)) +K \ \text{ as } \ \e \text{ is small enough};
\earray
\eeq
\beq{V-e}
\barray
\Ll^\e V^\e(x^\e(t),t)& = O(\e) (V(x^\e(t)) +1)
 + p \disp \sum_{i=1}^n [x_i^\e(t)]^p \bigg[ \lbar{b}_i + \frac{p}{2}\lbar{\sigma}^2_i -\sum_{j=1}^n \lbar{a}_{ij}x^\e_j(t) \bigg] \\
 &\leq O(\e) (V(x^\e(t)) +1) + p \disp \sum_{i=1}^n [x_i^\e(t)]^p \bigg[ \lbar{b}_i + \frac{p}{2}\lbar{\sigma}^2_i -\lbar{a}_{ii}x^\e_j(t) \bigg],
\earray
\eeq
where we used condition (A1).

Let $k_0 \in \nn$ be sufficiently large such that every component of $x^\e(0)$ is contained within the interval $\big( \dfrac{1}{k_0}, k_0\big)$. For each $k \geq k_0$, we define
\beq {ns}
\tau_k := \inf \bigg\{ t \in [0, \infty) : x^\e_i(t) \notin (\frac{1}{k}, k) \ \ \text{for some}\ i =1,2,\dots, n \bigg\}.
\eeq
Clearly, the sequence $\tau_k$, $k=1,2,\dots$ is monotonically increasing. Set $\tau_\infty := \lim_{k \to \infty} \tau_k$. We want to show that $\tau_\infty = \infty$ a.s. If this were false, there would exist some $T > 0 $ and $\tilde{\e} >0$ such that $P\{ \tau_\infty \leq T \} > \tilde{\e}$. Therefore, we can find some $k_1 \geq  k_0$ such that
\begin{equation}
P \{ \tau_k \leq T \} >\tilde{\e}, \ \ \text{for all} \ k \geq k_1.
\end{equation}
By \eqref{tV-e}, it can be verified that for any $(x,\al) \in \rr^n_+ \times \ZZ$,
\bea
\Ll^\e \tilde{V}^\e (x^\e(t),  t) \leq O(\e) \tilde{V}(x^\e(t)) + K
\eea
Using the generalized It\^{o}'s Lemma and taking the expectation on both sides, for any $k \geq k_1$, we have
\bea\ad
E^\e \tilde{V}^\e(x^\e(t \wedge \tau_k),  t \wedge \tau_k) - \tilde{V}^\e(x^\e(0), 0)
\leq \disp E^\e \int_0^{t \wedge \tau_k} O(\e) \tilde{V}(x^\e(s))ds + Kt.
\eea
Thus,
\bea
(1+O(\e)) E^\e \tilde{V} (x^\e(t \wedge \tau_k)) \leq \tilde{V}^\e(x^\e(0), \al^\e(0), 0)+ Kt+ \disp E^\e \int_0^{t \wedge \tau_k} O(\e) \tilde{V}(x^\e(s))ds.
\eea
When $\e$ is small enough, applying the generalized Gronwall's inequality, we obtain
\bea
 E^\e \tilde{V} (x^\e(t \wedge \tau_k)) \leq  \dfrac{\tilde{V}^\e(x^\e(0), \al^\e(0), 0)+ Kt}{1+ O(\e)}\ \  e^{\dfrac{O(\e)(t \wedge \tau_k)}{1+O(\e)} }.
\eea
Letting $t=T$, we have $E^\e \tilde{V} (x^\e(T \wedge \tau_k)) < \infty.$ On the other hand,
\bea \ad
E^\e \tilde{V} (x^\e(T \wedge \tau_k)) \geq E^\e [\tilde{V}(x^\e(\tau_k)) I_{\{ \tau_k \leq T\}}]\\
 \aad \quad > \tilde{\e} [(k^\ga -1 -\ga \log k)\wedge((1/k)^\ga -1+\ga \log k)] \to \infty,
\eea
as $k \to \infty$.
This is a contradiction so we must have $\lim_{k  \to \infty} \tau_k= \infty$ a.s.

By applying generalized It\^{o}'s Lemma to $e^t V^\e(x^\e(t), \al^\e(t), t)$ and taking the expectations of both sides, we have
\beq{eq2}
\barray
&\disp (1+O(\e))\bigg\{E^\e [e^{t \wedge \tau_k} \sum_{i=1}^n [x^\e_i(t \wedge \tau_k)]^p] -\sum_{i=1}^n [x^\e_i(0)]^p \bigg\}
\\
 \aad \quad = E^\e \int_0^{t \wedge \tau_k} e^s (V^\e(x^\e(s),s)+ \Ll^\e V^\e(x^\e(s), s) ds\\
 \aad \quad \leq E^\e \int_0^{t \wedge \tau_k}\Big[ pe^s \sum_{i=1}^n [x^\e_i(s)]^p \bigg( \frac{1+O(\e)}{p} +\lbar{b}_i + \frac{p}{2} \lbar{\sigma}_i^2 -\lbar{a}_{ii} x^\e_i(s) \bigg) + O(\e) e^s\Big] ds\\
\aad \quad  \leq E^\e \int_0^{t \wedge \tau_k} e^s K(\e) ds.
\earray
\eeq

By \eqref{eq2}, we have
$$E[e^{t \wedge \tau_k} \sum_{i=1}^n [x^\e_i(t \wedge \tau_k)]^p] -\sum_{i=1}^n [x^\e_i(0)]^p \leq E^\e \int_0^{t \wedge \tau_k} e^s K ds \leq K (e^t-1).
$$
Therefore, by virtue of Fatou's Lemma and  letting $k \to \infty$, we obtain
$$E \big[ \sum_{i=1}^n [x^\e_i(t)]^p  \big] \leq e^{-t} \sum_{i=1}^n [x^\e_i(0)]^p + K(1-e^{-t}) \leq K < \infty.$$
In view of the exponential dominance above, taking $\sup_{t \geq 0}$, we obtain the desired result.
 The next two parts of the theorem can be obtained similar to Section \ref{pros}. \qed

\subsection{Stability in Probability} \label{sub:stab}
Stability of dynamic systems with switching containing  randomly perturbed  processes has been done recently; see \cite{Bado}. In this study, our first goal is to establish the stability of \eqref{perturb} with small $\e$ via the stability of the averaged system \eqref{aver}.
We first recall the definition of stability for stochastic differential equations; see \cite{Khas3}.

\begin{defn}{\rm The equilibrium point $x=0$ of the system \eqref{aver} is said to be stable in probability, if for any $\e >0$ and any $\al \in \ZZ$,
$\lim_{y \ar 0}
P\{ \sup \limits_{t \geq 0}
|x^{y,\al}(t)|
> \e \}=0,$
where $x^{y,\al}(t)$ denotes the solution of \eqref{aver} with initial data $x(0)=y$ and $\al(0)=\al$.
}\end{defn}

Using  similar argument as \cite{Khas3}, we establish the following lemma.

\begin{lem}\label{stable}
Let $D \in \rr^n$ be a neighborhood of 0. Suppose that for each $i \in \ZZ$, there exists a non-negative function $V(\cdot,\al): D \mapsto \rr$ such that
\begin{itemize}
\item[{\rm(i)}] $V(\cdot,\al)$ is continuous in D and vanished only at $x=0$;
\item[{\rm(ii)}] $V(\cdot,\al)$ is twice continuously differentiable in $D \backslash \{0\}$ and
$ \Ll V(x,\al) \leq 0$,  $\forall  x \in D \backslash \{ 0\}.$
\end{itemize}
Then the equilibrium point $x=0$ is stable in probability.
\end{lem}

\begin{thm}
Assume that
\beq{ine}
(\lbar{r}_i - \lbar{a}_{ii})^2 +4 \lbar{a}_{ii} (\lbar{b}_i +\lbar{\sigma}_i^2) < 0 \text{, for all } i=1,2,\dots,n.
\eeq
Then under assumptions {\rm(A1)}, {\rm(A2)},  and \eqref{irre}, the equilibrium point $x=0$ is stable in probability for the averaged system \eqref{aver}.
\end{thm}

\para{Proof.}
We consider the Lyapunov function
\beq{lya}
V(x) = \disp \sum_{i=1}^n x_i - \log{(x_i+1)}.
\eeq
 It can be seen that $V(x)$ satisfies condition (i) of \lemref{stable}.

For  \eqref{aver}, we have
\beq{(A3)}
\Ll V(x) = \sum_{i=1}^n \frac{x_i^2}{x_i+1} \bigg( \lbar{r}_i - \sum_{j=1}^n \lbar{a}_{ij} x_j \bigg) + \frac{1}{2 } \sum_{i=1}^n \frac{1}{(x_i+1)^2} x_i^2 \lbar{\sigma}_i^2 .
 \eeq
 By  condition (A1), the property of solutions and the assumption, we have 
 \beq{neg} \barray
 \Ll V(x) & \leq \disp \sum_{i=1}^n \frac{x_i^2}{(x_i+1)^2} \Bigg\{ (x_i+1) (\lbar{r}_i - \lbar{a}_{ii} x_i)+ \frac{1}{2}\lbar{\sigma}_i^2 \Bigg\}\\
 & = \disp \sum_{i=1}^n \frac{x_i^2}{(1+x_i)^2} \Bigg[ -\lbar{a}_{ii}x_i^2 + (\lbar{r}_i -\lbar{a}_{ii})x_i+ (\lbar{b}_i + \lbar{\sigma}_i^2) \Bigg]
  < 0 \text{  for all  } x \neq 0.
\earray \eeq
Thus, by
\lemref{stable}, the equilibrium point $x=0$ of system \eqref{aver} is stable in probability. \qed

\begin{thm}\label{main} Under conditions {\rm(A1)}, {\rm(A2)}, \eqref{irre}, and \eqref{ine}, the equilibrium point $x=0$ is stable in probability for \eqref{perturb} for sufficiently small $\e$. 
\end{thm}

\para{Proof.}
With $V(x)$ defined by \eqref{lya}, $V^\e_1(x,  t)$ defined by \eqref{per-term}, $V^\e_2(x, t)$ defined by \eqref{per2} and their corresponding estimates, it is easy to see that $$V^\e(t)= V(x) +V^\e_1(x,  t)+V^\e_2(x, t)$$ satisfies condition (i) in  \lemref{stable}.
$V(x)$ is an increasing function and when $\e$ is small enough, by \thmref{Bound}, the process $x^\e(t)$ is stochastically bounded. Hence, $V(x^\e(t))$ is bounded for $\e$ is small enough.

Furthermore,
\beq{e-V}\barray
\Ll^\e V^\e(t) =\disp O(\e) \bigg[V(x^\e(t))+1\bigg] + \sum_{i=1}^n \frac{(x^\e_i (t))^2}{x^\e_i(t) +1} \big( \lbar{r}_i - \sum_{j=1}^n \lbar{a}_{ij} x_j^\e(t)\big)+\frac{1}{2}\sum_{i=1}^n \frac{(x^\e_i (t))^2}{(x^\e_i(t) +1)^2}\lbar{\sigma}_i^2.\earray
\eeq
By virtue of \eqref{neg}, $\Ll^\e V^\e(t) \leq 0$ for all $x^\e(t) \neq 0$ and $\e$ small enough. This verifies the theorem. \qed

\subsection{Extinction}\label{Extin}
In this section, we show if the averaged system \eqref{aver} is extinct, then the more complex switching system  \eqref{perturb}
is also extinct for sufficiently small $\e$.

\begin{defn}{\rm
The population is said to reach the extinction if $ \lim_{t \to \infty} |x(t)| =0$ a.s., i.e., $ \lim_{t\to \infty} \sum_{i=1}^n | x_i(t)| =0$ a.s.}
\end{defn}

\begin{thm}
Assume that
\beq{ine1}
\lbar{r}_i - \frac{1}{2} \lbar{\sigma}_i^2 \leq -c \text{, for all } i=1,2,\dots,n,
\eeq
where $c$ is a positive number.
Then under assumptions {\rm(A1)}, {\rm(A2)},  and \eqref{irre}, the population of the averaged system \eqref{aver} will become extinct exponentially a.s. for sufficiently small $\e$.
\end{thm}

 \para{Proof.}
 For each $i =1,2 \dots, n$, consider
 \beq{LY}
 V_i (x) = \log \left(x_i \right) .
 \eeq
 where $x_i$ is the $i$th component of $x$.
 Using the definition of the generator, we have
$$\Ll V_i(x(t)) = \lbar{r}_i -\sum_{j=1}^n \lbar{a}_{ij} x_j(t) - \frac{1}{2} \lbar{\sigma}_i^2. $$
Applying It\^{o}'s Lemma, we obtain
\bea
\log \left( x_i(t) \right) &=\log \left(x_i(0)\right) + \disp\int_0^t  \left(\lbar{r}_i -\sum_{j=1}^n \lbar{a}_{ij} x_j(s) - \frac{1}{2} \lbar{\sigma}_i^2\right) ds
+\int_0^t \lbar{\sigma}_i dw_i(s)\\
& \leq \log \left(x_i(0)\right) + t \left(\lbar{r}_i - \frac{1}{2} \lbar{\sigma}_i^2\right)
+ \lbar{\sigma}_i w_i(t).
\eea
$w_i(t)$ is a Brownian motion. Therefore, the strong law of large numbers for martingales implies that $\disp\lim_{t \to \infty}\dfrac{w_i(t)}{t} =0$ a.s. It follows by
$$\limsup_{t\to \infty} \frac{\log\left(x_i(t) \right)}{t} \leq \lbar{r}_i - \frac{1}{2} \lbar{\sigma}_i^2 \leq -c \quad a.s.$$
Thus, the sample Lyapunov exponent of the solution is negative, and the population will become extinct exponentially a.s. \qed

\begin{thm}\label{ext} Under conditions {\rm(A1)}, {\rm(A2)}, \eqref{irre}, and \eqref{ine1}, the population of the system  \eqref{perturb} will become extinct exponentially for sufficiently small $\e$.
\end{thm}

\para{Proof.}
With $V_i(x)$ defined by \eqref{LY}, $V^\e_{i,1}(x,  t)$ defined by \eqref{per-term}, $V^\e_{i,2}(x, t)$ defined by \eqref{per2} and their corresponding estimates, $V_i^\e (x) = V_i(x) + V^\e_{i,1}(x, t)+V^\e_{i,2}(x, t)$ satisfies the following properties:
\bea \ad V_i^\e(x^\e(t)) = V_i (x^\e(t) ) +O(\e) (V_i(x^\e(t)) +1 )\\
\ad \Ll^\e V_i^\e(x^\e(t)) = \disp O(\e) (V_i(x^\e(t))+1) + \lbar{r}_i -\frac{1}{2}\lbar{\sigma}_i^2 -\sum_{j=1}^n \lbar{a}_{ij} x^\e_j(t).
\eea
By the generalized It\^{o} Lemma,
\bea
\log \left(x_i^\e(t) \right) & = \disp \log \left( x_i^\e(0) \right) + \int_0^t \left[ O(\e) (V_i(x^\e(s))+1) + \lbar{r}_i -\frac{1}{2}\lbar{\sigma}_i^2 -\sum_{j=1}^n \lbar{a}_{ij} x^\e_j(s) \right] ds\\
&\qquad \  \disp + \int_0^t \sigma_i(\al^\e(s)) dw_i(s) \\
& \leq \disp  \log \left( x_i^\e(0) \right) + t \left[ O(\e) + \lbar{r}_i -\frac{1}{2}\lbar{\sigma}_i^2 \right] +O(\e) \int_0^t \log \left(x_i^\e(s) \right) ds\\
&\qquad \ \disp + \int_0^t \sigma_i(\al^\e(s)) dw_i(s).
\eea
Denote $M(t) = \int_0^t \sigma_i(\al^\e(s)) dw_i(s)$ and $M(t)$ is a martingale. Using the
quadratic variation of this martingale, we obtain that
$ t^{-1}
\langle M, M \rangle_t
= t^{-1} \int_0^t \sigma_i^2(\al^\e(s)) ds$ is bounded a.s.
The strong law of large numbers for martingales leads to
$\lim_{t \to \infty} {M(t)}/{t} =0$ a.s. (see \cite[Theorem 1.3.4]{Mao1}).
In addition, $\disp\lim_{t \to \infty} \frac{1}{t} \int_0^t \log \left( x_i^\e(s) \right) ds= \log (\lbar{x}_i(t))$ a.s, where $\lbar{x}_i(t)$ is the solution of \eqref{aver} (see \cite[Chapter 8]{YinZh}).Then $\disp\limsup_{t \to \infty} \frac{1}{t} \int_0^t \log \left( x_i^\e(s) \right) ds= \log (\lbar{x}_i(t))$ a.s.
Therefore,
\bea
\disp\limsup_{t \to \infty} \frac{\log \left( x^\e_i(t) \right)}{t}&\disp \leq O(\e)+ \lbar{r}_i -\frac{1}{2}\lbar{\sigma}_i^2 + \limsup_{t \to \infty} \frac{O(\e)}{t} \int_0^t   \log \left( x^\e_i(s) \right) ds\\
& \disp \leq O(\e)+ \lbar{r}_i -\frac{1}{2}\lbar{\sigma}_i^2 + O(\e) \log (\lbar{x}_i(t)) \ \hbox{ a.s.}
\eea
When $\e$ is small enough, under condition \eqref{ine1}, $\disp\limsup_{t \to \infty} \frac{\log \left( x^\e_i(t) \right)}{t} <0$ a.s. This results in the exponential extinction of the population. \qed

\subsection{Stochastic Permanence}\label{Perm}
We first recall the definition of stochastic permanence.

\begin{defn} {\rm The population system \eqref{aver} is said to be stochastically permanent if for any $\dl \in (0,1)$, there exist positive constants $H = H(\dl)$ and $K =K(\dl)$ such that
\beq{Pem}
\liminf_{t \to \infty} P \{ |x(t)| \geq H \} \geq 1-\dl, \quad \liminf_{t \to \infty} P \{ |x(t)| \leq K \} \geq 1-\dl,
\eeq
where $x(t)$ is the solution of the population system \eqref{aver} with any initial condition $x(0) \in \rr ^n_+$.
}\end{defn}

 \begin{lem}\label{Per} Assume that {\rm(A1)}, {\rm(A2)}, and \eqref{irre} hold. Then  population  system \eqref{aver} is stochastically permanent when $\lbar{b}_i > 0$ for $i =1,2,\dots,n$.
 \end{lem}

 \para{Proof.}
 To obtain the stochastic permanence, we need to prove two inequalities in \eqref{Pem} and the first part is followed by \thmref{Bound}.
 Before working on the second part, we first set the notation: $\wdt{r} := \max \lbar{r}_i$, $\hat{r}=\min \lbar{r}_i$, $\hat{b}=\min \lbar{b}_i$, $\wdt{a} = \max \lbar{a}_{ij}$, $\wdt{\sg} = \lbar{\sg}_i$.
  We begin to work with some estimates for the averaged system \eqref{aver}, where $x(t)$ is the solution.
 Let $\theta$ be a positive constant such that
 $\theta \wdt{\sg}^2 < 2 \hat{b},$
 and $\kappa > 0$ satisfying $0 < \dfrac{2 \kappa}{\theta} < 2 \hat{b} - \theta \wdt{\sg}^2$.
Consider
\bea\ad V(x) = \disp\sum_{i=1}^n x_i,\
U(x) = \dfrac{1}{V(x)}, \ \hbox{ and }\
 J(x) = e^{\kappa t} \big(1+ U(x) \big)^\theta .\eea
By applying It\^{o}'s Lemma, we have
\bea dU(x(t))\ad = \bigg[-U^2(x(t)) \sum_{i=1}^n x_i (t)(\lbar{r}_i - \sum_{j=1}^n \lbar{a}_{ij}x_j(t)) + U^3(x(t)) \sum_{i=1}^n \lbar{\sg}^2_i x^2_i(t) \bigg] dt \\
\aad \quad - U^2(x(t)) \sum_{i=1}^n \lbar{\sg}_i x_i(t) d w_i(t).\eea
Note that
\beq{J}
\barray
 dJ(x(t))
 \ad  = \disp \theta e^{\kappa t}(1+U(x(t)))^{\theta -2} \bigg\{ \bigg[\frac{\kappa}{\theta} (1+U(x(t)))^2 -(1+U(x(t))) U^2(x(t)) \sum_{i=1}^n x_i(t)\\
 \aad \quad \times (\lbar{r}_i -\sum_{j=1}^n \lbar{a}_{ij} x_j(t))   + U^3(x(t)) \sum_{i=1}^n \lbar{\sg}^2_i x_i^2(t) + \frac{\theta +1}{2} U^4(x(t)) \sum_{i=1}^n \lbar{\sg}^2_i x_i^2(t)\bigg] dt\\
 \aad \quad - (1+U(x(t))) U^2(x(t)) \sum_{i=1}^n \lbar{\sg}_i x_i(t) dw_i(t) \bigg\}.
\earray
\eeq

Therefore,
\beq{eq5}
\barray
\Ll J(x) \ad = \disp \theta e^{\kappa t} (1+U(x))^{\theta -2} \bigg[\frac{\kappa}{\theta} (1+U(x))^2 -(1+U(x)) U^2(x) \sum_{i=1}^n x_i (\lbar{r}_i -\sum_{j=1}^n \lbar{a}_{ij} x_j) \\
\aad \ + U^3(x) \sum_{i=1}^n \lbar{\sg}^2_i x_i^2 + \frac{\theta +1}{2} U^4(x) \sum_{i=1}^n \lbar{\sg}^2_i x^2_i \bigg].
\earray
\eeq
We have
\beq{term1}
\barray
\ad \!\!\!\! -(1+U(x(t))) U^2(x(t)) \sum_{i=1}^n x_i (t)(\lbar{r}_i -\sum_{j=1}^n \lbar{a}_{ij} x_j(t))\\
  \aad\leq  \disp -(1+U(x(t))) U^2(x(t)) \sum_{i=1}^n x_i (t)(\lbar{r}_i - \wdt{a} \sum_{j=1}^n x_j(t))\\
 \aad =\disp -(1+U(x(t))) U^2(x(t)) \bigg[\sum_{i=1}^n x_i(t) (\lbar{b}_i + \frac{1}{2} \lbar{\sg}_i^2) - \wdt{a} (\sum_{i=1}^n x_i(t)) (\sum_{j=1}^n x_j(t))\bigg]\\
 \aad \leq \disp -(1+U(x(t))) U^2(x(t)) \sum_{i=1}^n x_i (t)\hat{b} -(1+U(x(t))) U^2(x(t))\sum_{i=1}^n \frac{1}{2}x_i (t)\lbar{\sg}^2_i\\
  \aad \qquad + (1+U(x(t))) \wdt{a} \\
\aad \leq \disp -\hat{b} (1+U(x(t)))U(x(t)) - U^3(x(t)) \sum_{i=1}^n x_i (t)\frac{\lbar{\sg}^2_i}{2} +(1+U(x(t))) \wdt{a}\\
\aad \leq \disp -\hat{b} (1+U(x(t)))U(x(t))  -\frac{1}{2}U^4(x(t)) \sum_{i=1}^n x^2_i(t) \lbar{\sg}^2_i +(1+U(x(t))) \wdt{a}.
\earray
\eeq
\beq{term2}
 U^3(x(t)) \disp \sum_{i=1}^n \lbar{\sg}^2_i x^2_i(t) \leq U^3(x(t)) \wdt{\sg}^2 \disp \sum_{i=1}^n x^2_i(t)
 \leq U(x(t)) \wdt{\sg}^2 \dfrac{\disp \sum_{i=1}^n x^2_i(t)}{(\disp \sum_{i=1}^n x_i(t))^2}
 \leq \wdt{\sg}^2 U(x(t)).
\eeq
Thus,
\beq {eq6}
\barray
 \Ll J(x(t))
 \ad\leq \disp \theta  e^{\kappa t}(1+U(x(t)))^{\theta -2}\bigg[\frac{\kappa}{\theta} (1+U(x(t)))^2 -\hat{b} (1+U(x(t)))U(x(t))  \\
\aad  -{U^4(x(t))\over 2} \sum_{i=1}^n x^2_i(t) \lbar{\sg}^2_i+(1+U(x(t))) \wdt{a} + \wdt{\sg}^2 U(x(t))\! +\! \frac{\theta +1}{2} U^4(x(t)) \sum_{i=1}^n \lbar{\sg}^2_i x^2_i (t)  \bigg]\\
  \ad\leq \disp \theta e^{\kappa t} (1+U(x(t)))^{\theta -2}\bigg[\frac{\kappa}{\theta} + \frac{2\kappa}{\theta} U(x(t)) + \frac{\kappa}{\theta}  U^2(x(t)) -\hat{b} U(x(t))-\hat{b} U^2(x(t))) \\
  \aad + \wdt{a}+ \wdt{a} U(x(t)) +\wdt{\sg}^2 U(x(t)) + \frac{\theta}{2} \wdt{\sg}^2 U^2(x(t)) \bigg]
  \leq K \theta  e^{\kappa t}.
\earray
\eeq
where $K$ is a positive constant  depending on $\kappa$, $\theta$, and coefficients of the system. (This inequality is resulted from the choice of $\theta$ and $\kappa$.)

Integrating and taking expectations on both sides of \eqref{J}, we have:
$
E[J(x(t))] -J(x(0)) \leq K \int_0^t \theta  e^{\kappa t} ds,
$
i.e.,
$$
E[(1+U(x(t)))^\theta] \leq e^{-\kappa t} (1+U(x(0)))^\theta +\frac{K \theta}{\kappa}.
$$
Note that for $x \in \rr^n_+$, $\disp (\sum_{i=1}^n x_i)^\theta \leq n^\theta |x|^\theta$. For any given $\dl \in (0,1)$, choose $H >0$ such that $\dfrac{H^\theta n^\theta K \theta}{\kappa} \leq \dl.$ By Tchebychev's inequality, we obtain
$$
\barray
P(|x(t)| <H) &\leq P \bigg( U^\theta(x(t)) > \dfrac{1}{H^\theta n^\theta}  \bigg)\\
& \leq H^\theta n^\theta E[U(x(t))^\theta] \\
& \leq H^\theta n^\theta E[(1+U(x(t)))^\theta] \\
& \leq H^\theta n^\theta \bigg[e^{-\kappa t} (1+U(x(0)))^\theta +\frac{K \theta}{\kappa} \bigg].
\earray
$$
This implies that $\disp \limsup_{t \to \infty} P(|x(t)| < H) \leq \frac{H^\theta n^\theta K \theta}{\kappa} \leq \dl$, i.e. $\disp \liminf_{t \to \infty} P(|x(t)| \geq  H) \geq 1-\dl$. This completes the proof. \qed

\begin{thm}
Under  conditions {\rm(A1)}, {\rm(A2)}, and \eqref{irre}, system \eqref{perturb} is stochastically permanent when $\lbar{b}_i >0$ for each $i=1,2,\dots,n$ and sufficiently small $\e$.
\end{thm}

\para{Proof.}
We apply the definition of the generator \eqref{gen} and obtain the following for the perturbed system \eqref{perturb}, for each $\alpha$,
\bea
\Ll^\e J(x)\ad = \sum_{i=1}^n J_{x_i}(x) \xi_i(x, \al) + \frac{1}{2} \sum_{i=1}^n J_{x_i x_i} (x) s^2_i (x, \al).
\eea
Similar to \eqref{per-term} and \eqref{per2}, we define
\bea
  J^\e_1(x,t)\ad =E^\e_t \int_t^\infty e^{t-u} \disp\sum_{i=1}^n J_{x_i}(x) [\xi_i(x,\al^\e(u)) -\lbar{\xi}_i(x)] du,\\\\
  J^\e_2(x,t)\ad =E^\e_t \int_t^\infty e^{t-u} \disp\frac{1}{2}\sum_{i=1}^n J_{x_i x_i}(x) [s_i^2(x, \al^\e(u)) -\lbar{\lambda}^2_i(x)] du.
\eea
Then
\bea \ad J^\e_1(x^\e(t),t) = O(\e) [J(x^\e(t)) +1], \\
\ad J^\e_2(x^\e(t),t) = O(\e) [J(x^\e(t)) +1],\\
\ad \Ll^\e J^\e_1(x^\e(t), t) = O(\e)(J(x^\e(t),t)+1)\\
\aad \qquad - \disp \sum_{i=1}^n J_{x_i} (x^\e(t))[\xi_i(x^\e(t),\al^\e(t)) -\lbar{\xi}_i (x^\e(t))],\\
\ad \Ll^\e J^\e_2(x^\e(t), t) = O(\e)(J(x^\e(t),t)+1)\\
\aad \qquad - \disp \sum_{i=1}^n J_{x_i x_i} (x^\e(t))[s^2_i(x^\e(t),\al^\e(t)) -\lbar{\lambda}^2_i (x^\e(t))].
\eea
Define $$J^\e (x,t) = J(x) + J^\e_1 (x,t)+J^\e_2 (x, t).$$ The functions satisfy
\beq{J-e}
\barray
J^\e(x^\e(t),t)\ad = J(x^\e(t))+ O(\e) (J(x^\e(t)) +1)\\
\Ll^\e J^\e(x^\e(t),t)\ad = O(\e) (J(x^\e(t)) +1)
 + \theta e^{\kappa t} (1+U(x^\e(t)))^{\theta -2} \bigg[\frac{\kappa}{\theta} (1+U(x^\e(t)))^2 \\
 \aad\quad - (1+U(x^\e(t))) U^2(x^\e(t)) \sum_{i=1}^n x^\e_i(t) (\lbar{r}_i -\sum_{j=1}^n \lbar{a}_{ij} x^\e_j(t))\\
 \aad \quad + U^3(x^\e(t)) \sum_{i=1}^n \lbar{\sg}^2_i (x^\e_i(t))^2 + \frac{\theta +1}{2} U^4(x^\e(t)) \sum_{i=1}^n \lbar{\sg}^2_i (x^\e_i(t))^2\bigg]\\
 \ad = O(\e) (J(x^\e(t)) +1) + \theta e^{\kappa t} K.
\earray
\eeq
Integrating  on both sides of \eqref{J-e} and taking expectation, we have
$$
E^\e_t[J^\e(x^\e(t), t)]-J^\e(x^\e(0), 0) \leq O(\e)E^\e_t [J(x^\e(t))] + O(\e) + \frac{\theta K}{\kappa} e^{\kappa t}.
$$
Denote $J^\e_0 = J^\e(x^\e(0), 0)$. Then
\bea
(1-O(\e)) E^\e_t[J(x^\e(t))] \leq O(\e) + J^\e_0 +\dfrac{\theta K}{\kappa} e^{\kappa t},
\eea
i.e., when $\e >0$ is small enough, $1 - O(\e) >0$ and
$$
E^\e_t[(1+U(x^\e(t)))^\theta] \leq \frac{O(\e) +J^\e_0}{1-O(\e)} e^{-\kappa t} + \frac{\theta K}{\kappa}.
$$
For any given $\dl \in (0, 1)$, choose $H >0$ such that $\disp \frac{H^\theta n^\theta K \theta}{\kappa} \leq \dl$, by using Tchebychev's inequality, we can obtain the first inequality in \eqref{Pem}. The second inequality in \eqref{Pem} is obtained from \thmref{Bound}. This completes the proof. \qed

\section{Concluding Remarks}\label{Rem}
The paper has been devoted to competitive Lotka-Volterra model.
\begin{itemize}
\item We formulate the ecosystems as a hybrid systems involve both continuous states and discrete events in which the discrete events take values in a countable state space.
We demonstrated such properties as
 existence and uniqueness of solution, stochastic boundedness, sample path continuity for the models.

\item A main effort is placed on reduction of complexity by introducing a small parameter into the system. This leads to a  two-time-scale formulation.
Although the two-time-scale system has complex structures, it is shown that there is an associated averaged or reduced system.
\item
Using the averaged system, we prove that extinction and permanence for the Lotka-Volterra ecosystems with a two-time-scale Markov chain
   (of the complex original system)
    by perturbed Lyapunov function methods when the $\e$ is small enough.

\item
 A number of questions deserve further consideration.
  \begin{itemize}
  \item To begin, instead of the current formulation, we may consider the Markov chain involves both fast and slow motions with more complex structure.
 For example, two-time-scale Markov chains that are nearly decomposable were considered in \cite{YinHZ}. Such setups may be adopted to the ecosystems.

  \item
 Other related systems such as mutualism systems can also be formulated and studied.
 Moreover, one may consider populations suffering sudden environmental shock (e.g., earthquakes, hurricanes, tornadoes, etc.), leading to the consideration regime-switching jump diffusion systems. Designing feedback controls so as to achieve permanence and extinction etc.
is another area of future study.
\item There is a growing interest to study the associate harvesting problems \cite{TY15}.
To study the harvesting strategies with systems proposed in this paper has not been done to date and is a worthwhile direction.
\end{itemize}
\end{itemize}

 \appendix
 \section{Proofs of Technical Results}

 \para{Proof of Theorem \ref{supr}.}
 The proof consists of two parts. In the first part, we show that there is a unique global solution, and in the second part, we show the solution lives in $\rr^n_+$.

Step 1:  For any $\iota \in \ZZ$, in view of \cite[Theorem 2.1]{MMR}
there is a unique strong solution for the following diffusion
\beq{neq}
dx(t) = \Xi (x(t),\iota)dt + S(x(t),\iota) dw(t), \ x(0)= x_0 \in \rr^n_+.
\eeq
The rest of the proof of this part is similar to that of \cite[Theorem 3.1]{NY16}, so we will be brief.
 For any  stopping time $\tau$ and an $\F_\tau$-measurable $\rr^n$-valued random variable $x(\rho)$, there exists a strong  solution to \eqref{neq} in $[\rho,\infty)$; see \cite[Remark 3.10]{MaoY}.  We proceed to construct the solution  with any initial data $(x_0,i_0)  \in \rr_+^n \times \ZZ$ by the interlacing procedure \cite[Chapter 5]{Apple}. Denote by $\tilde{x}^{(0)}(t), t \geq 0$  the solution to
$$ d\tilde{x}^{(0)}(t) =\Xi (\tilde{x}^{(0)}(t),i_0)dt + S(\tilde{x}^{(0)}(t),i_0) dw(t), \
\tilde x^{(0)} (0)= x_0.$$
Set
$\rho_1 = \inf \{ t >0: \int_0^t \int_\rr h(i_0,z)\p(ds,dz) \neq 0\}$,
$i_1 = i_0 + \int_0^{\rho_1} \int_\rr h(i_0,z) \p(ds,dz),
$
and let $\tilde{x}^{(1)}(t), t \leq \rho_1$ be the solution  to
$$d\tilde{x}^{(1)}(t) =\Xi (\tilde{x}^{(0)}(t),i_1)dt + S(\tilde{x}^{(0)}(t),i_1) dw(t),$$
with  initial data $\tilde{x}^{(1)}(\rho_1)=\tilde{x}^{(0)}(\rho_1)$.
Continuing this procedure, let $\disp \rho_\infty = \lim_{k \to \infty} \rho_k$ and set
$x(t) = \tilde{x}^{(k)}(t)$, $\al(t) =i_k$, if  $\rho_k \leq t < \rho_{k+1}.$
Then
$$
\begin{cases}
\disp
x(t \wedge \rho_k)= x_0 + \int_0^{t \wedge \rho_k}\Xi (x(s),\al(s))ds + S(x(s),\al(s)) dw(s), \\ \disp
\al(t \wedge \rho_k) =i_0+ \int_0^{t \wedge \rho_k} \int_\rr h(\al(s-),z) \p(ds,dz),
\end{cases}
$$
where $\p(ds,dz)$ is a Poisson random measure as defined in \cite[p. 29]{YinZ10} with modification to countable state space; see also \cite{NY16}.
 To verify that $x(t)$ is a global solution,
 we claim that $\rho_\infty = \infty$.
 In fact, it can be shown as in \cite{NY16},
  for any $T >0$,$
P( \rho_k \leq T )
\leq \disp \sum_{l=k}^\infty e^{-MT}\frac{(MT)^l}{l!}.
$
Thus
$P( \rho_k \leq T )  \to 0$ as $k \to \infty$ so $\rho_\infty =\infty$  a.s.
The uniqueness of $x(t)$ follows from the uniqueness of $\tilde{x}^{(k)}(t)$ on $[\rho_k, \rho_{k+1})$.
Thus, we have shown that there is a unique global solution to
$dx(t) = \Xi (x(t),\alpha(t))dt + S(x(t),\alpha(t)) dw(t)$
with arbitrary initial data $(x_0,i_0)$.

 Step 2: Show the solution  $x(t)$ obtained in Step 1 above remains in $\rr^n_+$.
The proof is similar to
 \cite[Theorem 2.1]{YinZ} although the switching set is now countable. Let $k_0 \in \nn $ be sufficiently large such that every component of $x(0)$ is contained in  $(\frac{1}{k_0}, k_0)$. For each $k \geq k_0$,  define
\beq{time}
\zeta_k := \inf \bigg\{ t \in [0, \zeta): x_i(t) \notin (\frac{1}{k},k) \text{ for some } i=1,2,\dots,n\bigg\}.
\eeq
The sequence ${\zeta_k, k=1,2,\dots}$ is monotonically so there is a limit $\zeta_{\infty}: = \lim_{k \rightarrow \infty} \zeta_k$ with $\zeta_\infty \leq \zeta$. We are to show $\zeta_{\infty} = \infty$ a.s. For suppose not, there would exist some $T>0$ and $\e >0$ such that $P\{ \zeta_\infty \leq T \} > \e$. Therefore, we can find some $k_1 \geq k_0$ such that
\beq{contr}
P \{ \zeta_k \leq T \} >\e, \quad \hbox{for all } \quad  k \geq k_1.
\eeq
Now, we consider the following Lyapunov function $V(x, \al) =V(x)$ independent of $\al$ given by
$V(x)
= \sum_{i=1}^n [x_i^\gamma -1 - \gamma \log{x_i}]$ for  $x \in \rr^n_+$ and  $0 < \gamma <1.$
Detailed calculation shows that for all $x \in \rr^n_+ $, $V(x) \geq 0$ and
$
\Ll V(x)
\leq K < \infty,$
where in the above, we used condition (A1). In view of  It\^{o}'s Lemma \cite{Sko}, for any $k \geq k_1$,
$$V(x(T \wedge \zeta_k)) - V(x(0)) = \int_0^{T \wedge \zeta_k} \Ll V(x(s))ds+ \sum_{i=1}^n \int_0^{T \wedge \zeta_k} \gamma \sigma_i(\al(s))(x_i^\gamma(s)-1)dw_i(s).$$
By virtue of Dynkin's formula and the bound $\Ll v(x) \le K$,
$KT + V(x(0)) \geq E[V(x(T \wedge \zeta_k))] \geq E[V(x(\zeta_k)) I_{\{\zeta_k \leq T\}}].$
By the definitions of $\zeta_k$ and V, we have
$V(x(\zeta_k) ) \geq (k^\gamma -1 -\gamma \log k) \wedge(\frac{1}{k^\gamma} -1 +\gamma \log k),$
and hence, it follows from \eqref{contr} that
$$\barray
KT + V(x(0)) &\geq [(k^\gamma -1 -\gamma \log k) \wedge(\frac{1}{k^\gamma} -1 +\gamma \log k)] P\{\zeta \leq T\}
\ar \infty, \hbox{ as } k \ar \infty.
\earray$$
This is a contradiction, so we must have $\lim_{k \ar \infty} \zeta_k = \infty$ a.s., so $\zeta = \infty$ a.s. Thus,  the solution of \eqref{neq} remains in $\rr_+^n$ almost surely.
\qed

\para{Proof of Theorem \ref{moment}.}
Let $k_0 \in \nn $ be sufficiently large such that every component of $x(0)$ is contained in the interval $(\frac{1}{k_0}, k_0)$. For each $k \geq k_0$, we define
$\tau_k := \inf \{ t \in [0, \infty): x_i(t) \notin (\frac{1}{k},k)$ for some  $i=1,2,\dots,n\}$.
Similar to the proof in Step 2 of Theorem \ref{supr}, we can show that
$\lim_{k \ar \infty} \tau_k = \infty$ a.s.

Consider $V(x) = \disp \sum_{i=1}^n x_i^p$. Then it follows that for $x \in \rr^n_+$, we have
\beq{mom}\barray
\Ll V(x)\ad = p \disp\sui
\sum_{i=1}^n x_i^p \Big[ b_i(\iota) + \frac{p}{2} \sigma_i^2(\iota) -\sum_{j=1}^n a_{ij}(\iota) x_j \Big]\oal\\
\ad
 \leq p \sui\sum_{i=1}^n x_i^p \Big[ b_i(\iota) + \frac{p}{2} \sigma_i^2(\iota) - a_{ii}(\iota) x_i \Big]\oal ,\earray
\eeq
where in the last step, we used condition (A1).
By applying generalized It\^{o}'s Lemma \cite{Sko} to $e^t V(x(t))$, we have
\bea \ad e^{t \wedge \tau_k} \sum_{i=1}^n x_i^p(t \wedge \tau_k)- \sum_{i=1}^n x_i^p(0) \\
 \aad \ = \int_0^{t \wedge \tau_k} e^s (V(x(s))+\Ll V(x(s))) ds +p \sum_{i=1}^n \int_0^{t\wedge \tau_k} e^s x_i^{p-1}(s) \sigma_i(\al(s)) dw_i(s) ,\eea
where $\tau_k$ is the stopping time defined at the beginning of the proof. Thus taking expectations on both sides and using the assumption (H1), we obtain from \eqref{mom} that
\beq{eq1}
 E[e^{t \wedge \tau_k} \sum_{i=1}^n x_i^p(t \wedge \tau_k)]- \sum_{i=1}^n x_i^p(0)
 = E\int\limits_0^{t \wedge \tau_k} e^s (V(x(s))+\Ll V(x(s))) ds
 \leq E \int\limits_0^{t \wedge \tau_k} e^s Kds .
 \eeq
 By \eqref{eq1}, we have
$$E[e^{t \wedge \tau_k} \sum_{i=1}^n x_i^p(t \wedge \tau_k)]- \sum_{i=1}^n x_i^p(0) \leq E \int\limits_0^{t \wedge \tau_k} e^s Kds \leq K(e^t-1).$$
 Therefore, by virtue of Fatou's Lemma  and letting $k \ar \infty$, we obtain that
$$E\bigg[\sum_{i=1}^n x_i^p(t) \bigg] \leq e^{-t} \sum_{i=1}^n x_i^p(0) + K(1 - e^{-t})  \leq K < \infty.$$
In view of the exponential dominance above, taking $\sup_{t \geq 0}$, we obtain the desired result. \qed

\para{Proof of Theorem \ref{continu}.}
For any $0 \leq \tilde{t} \leq t$, we have
$$x_i(t)-x_i(\tilde{t}) = \int_{\tilde{t}}^{t} {\xi_i(x(r), \al(r))dr} + \int_{\tilde{t}}^{t} s_i(x(r),\al(r))dr, $$
and hence
\beq{dif}
|x_i(t) -x_i(\tilde{t})|^4 \leq 8 \bigg| \int_{\tilde{t}}^t \xi_i(x(r), \al(r))dr  \bigg|^4 + 8 \bigg| \int_{\tilde{t}}^t s_i(x(r),\al(r))dr\bigg|^4.
\eeq
Detailed computations in \thmref{moment} and H\"older's inequality lead to
\beq{ter1}
E \bigg| \int_{\tilde{t}}^t \xi_i(x(r), \al(r))dr  \bigg|^4 \leq (t-\tilde{t})^3 E \int_{\tilde{t}}^t |\xi_i(x(r),
\al(r))|^4 dr \leq K |t-\tilde{t}|^4.
\eeq
Meanwhile, we can
 show that
$E \bigg| \int_{\tilde{t}}^t s_i(x(r),\al(r))dr\bigg|^4 \leq K |t-\tilde{t}|^2.$
Thus
\beq{con}
E[|x(t)-x(\tilde{t})|^4] \leq K |t-\tilde{t}|^2.
\eeq
The desired result then follows from the well-known Kolmogorov continuity criterion. \qed

\para{Proof of Theorem \ref{O-e}.} The proof here is similar to
\cite{Bado}.
Direct calculations leads to
$$ \barray\disp
E \Bigg[\int_0^\infty e^{-t} (I_{\{\al^\e(t) =\al\}}-\nu_\al) dt\Bigg]^2 &
 \disp= \Big[\int_0^\infty \int_0^t e^{-t-s} +
  + \int_0^\infty \int_0^s\Big]
  O(\e + e^{-\kappa_0 (t-s)/\e})dsdt.
 \earray $$
Furthermore, $O(\e) \int_0^\infty \int_0^t e^{-t-s} dsdt = O(\e)$.
In addition, for some $K>0$,
$$\barray
\disp\int_0^\infty \int_0^t e^{-t-s} O(e^{-\kappa_0 (t-s)/\e})dsdt &\disp \leq K \int_0^\infty \int_0^t e^{-t (\kappa_0+\e)/\e} e^{s(\kappa_0-\e)/\e}dsdt
\leq K \dfrac{\e}{2 (\kappa_0 -\e)} = O(\e).
\earray$$
Thus, $\int_0^\infty \int_0^t e^{-t-s} O(\e + e^{-\kappa_0 (t-s)/\e})dsdt = O(\e)$.
 Likewise, by symmetry, we also have $\int_0^\infty \int_0^s e^{-t-s} O(\e + e^{-\kappa_0 (s-t)/\e})dtds =O(\e)$. The proof is complete. \qed

\end{document}